\documentclass[a4paper,12pt]{amsart}

\usepackage{amsmath}
\usepackage{amssymb}
\usepackage{mathrsfs}
\usepackage{ifthen}
\usepackage{graphicx}
\usepackage[T1]{fontenc} 

\nonstopmode \numberwithin{equation}{section}
\newtheorem{thm}{Theorem}[section]
\newtheorem{cor}[equation]{Corollary}
\newtheorem{lem}[equation]{Lemma}

\theoremstyle{definition}
\newtheorem{defn}{Definition}[section]

\newtheorem{prob}[equation]{Problem}

\newcounter{minutes}
\divide\time by 60
\newcounter{hours}
\multiply\time by 60
\addtocounter{minutes}{-\time}

\newcounter {own}
\def\theown {\thesection       .\arabic{own}}

{\qed\bigskip}

\newcounter{alphabet}

\def\be{\begin{equation}}
\def\ee{\end{equation}}

\newcommand{\bee}{\begin{enumerate}}
\newcommand{\eee}{\end{enumerate}}

\newcommand{\blem}{\begin{lem}}
\newcommand{\elem}{\end{lem}}
\newcommand{\bthm}{\begin{thm}}
\newcommand{\ethm}{\end{thm}}
\newcommand{\bcor}{\begin{cor}}
\newcommand{\ecor}{\end{cor}}
\newcommand{\beg}{\begin{examp}}
\newcommand{\eeg}{\end{examp}}
\newcommand{\begs}{\begin{examples}}
\newcommand{\eegs}{\end{examples}}
\newcommand{\bdefe}{\begin{defin}}
\newcommand{\edefe}{\end{defin}}
\newcommand{\bprob}{\begin{prob}}
\newcommand{\eprob}{\end{prob}}
\newcommand{\bei}{\begin{itemize}}
\newcommand{\eei}{\end{itemize}}

\begin{document}

\title{Proof of The Generalized Zalcman Conjecture for Initial Coefficients of Univalent Functions}

\author{Vasudevarao Allu}
\address{Vasudevarao Allu,
Discipline of Mathematics,
School of Basic Sciences,
Indian Institute of Technology  Bhubaneswar,
Argul, Bhubaneswar, PIN-752050, Odisha (State),  India.}
\email{avrao@iitbbs.ac.in}

\author{Abhishek Pandey}
\address{Abhishek Pandey,
Discipline of Mathematics,
School of Basic Sciences,
Indian Institute of Technology  Bhubaneswar,
Argul, Bhubaneswar, PIN-752050, Odisha (State),  India.}
\email{ap57@iitbbs.ac.in}

\subjclass[2010]{Primary 30C45, 30C50}
\keywords{Analytic, univalent,  starlike, convex, functions, coefficients, holomorphic motions,  variational method, Zalcman conjecture, }

\def\thefootnote{}
\footnotetext{ {\tiny File:~\jobname.tex,
printed: \number\year-\number\month-\number\day,
          \thehours.\ifnum\theminutes<10{0}\fi\theminutes }
} \makeatletter\def\thefootnote{\@arabic\c@footnote}\makeatother

\thanks{}

\maketitle
\pagestyle{myheadings}
\markboth{Vasudevarao  Allu and Abhishek Pandey}{Proof of The Generalized Zalcman Conjecture for Initial Coefficients}

	\begin{abstract}
		Let $\mathcal{S}$ denote the class of analytic and univalent ({\it i.e.}, one-to-one) functions $	f(z)= z+\sum_{n=2}^{\infty}a_n z^n$ in the unit disk $\mathbb{D}=\{z\in \mathbb{C}:|z|<1\}$. For $f\in \mathcal{S}$,  Ma proposed the generalized Zalcman conjecture that  $$|a_{n}a_{m}-a_{n+m-1}|\le (n-1)(m-1),\,\,\,\mbox{ for } n\ge2,\, m\ge 2,$$
		with equality  only for the Koebe function $k(z) = z/(1 - z)^2$ and its rotations.  In this paper using the properties of holomorphic motion and  Krushkal's Surgery Lemma  \cite{Krushkal-1995}, we prove the generalized Zalcman conjecture  when $n=2$, $m=3$ and $n=2$, $m=4$. 
	
	\end{abstract}

	\section{Introduction and Preleminiries}
	Let $\mathcal{H}$ denote the class of analytic functions in the unit disk $\mathbb{D}:=\{z\in\mathbb{C}:\, |z|<1\}$.  Let $\mathcal{A}$ be the class of functions 
	$f\in \mathcal{H}$ such that $f(0)=0$ and $f'(0)=1$, and denote by $\mathcal{S}$ 
	 the class of functions $f\in\mathcal{A}$ which are univalent ({\it i.e.}, one-to-one) in $\mathbb{D}$. 
	Thus $f\in\mathcal{S}$ has the following representation
	\begin{equation}\label{p4_i001}
		f(z)= z+\sum_{n=2}^{\infty}a_n z^n.
	\end{equation}
	
	Closely related to the class $\mathcal{S}$ is the class $\sum$ of functions which are analytic and univalent in the domain $\Delta=\{z: |z|>1\}$ exterior to $\mathbb{D}$, except for simple pole at infinity with residue $1$ and $g(z)\ne 0$ in $\Delta$. For $f\in \mathcal{S}$ given by (\ref{p4_i001}), let $F\in \Sigma$ be given by
	\begin{equation}\label{class sigma}
		F(z)=\frac{1}{f\left(\dfrac{1}{z}\right)}=z+b_{0}+\sum_{n=1}^{\infty}b_{n}z^{-n}.
	\end{equation}
	
	A simple computation in (\ref{class sigma}) shows that
	\begin{eqnarray}\label{S to sigma}
		& a_{2}=& -b_{0}\\[2 mm] \nonumber
		& a_{3}=& -b_{1}+b_{0}^2\\[2 mm] \nonumber
		&a_{4}=& -b_{2}+2b_{1}b_{0}-b_{0}^3\\[2 mm] \nonumber
		&a_{5}=& -b_{3}+2b_{2}b_{0}+b_{1}^2-3b_{1}b_{0}^2+b_{0}^4, 
	\end{eqnarray}
	
	and so

	\begin{eqnarray}\label{equ-p4-1}
		 a_{2}a_{3}-a_{4}&=&-b_{0}b_{1}+b_{2}\\[2 mm] 
		 a_{2}a_{4}-a_{5}&=&-b_0b_2+b_1b_0^2+b_{3}-b_{1}^2.
	\end{eqnarray}
	\vspace{1mm}

	A set $\Omega \subseteq \mathbb{C}$ is said to be starlike with respect to a point $z_{0}\in \Omega$ if the line segment joining $z_{0}$ to every other point $z\in \Omega$ lies entirely in $\Omega$. A set $\Omega\subseteq\mathbb{C}$ is called convex if the line segment joining any two points of $\Omega$ lies entirely in $\Omega$. A function $f\in\mathcal{A}$ is called starlike (convex respectively) if $f(\mathbb{D})$ is starlike with respect to the origin (convex respectively). Let $\mathcal{S}^*$ and $\mathcal{C}$  denote the classes of starlike and convex functions in $\mathcal{S}$ respectively. Then it is well-known that a function $f\in\mathcal{A}$ belongs to $\mathcal{S}^*$ if,  and only if, ${\rm Re\,}\left(zf'(z)/f(z)\right)>0$ for $z\in\mathbb{D}$. A function $f\in \mathcal{A}$ is called close-to-convex if there exists a real number $\theta$ and a function $g \in \mathcal{S}^*$ such that  $\Re (e^{i\theta}zf'(z)/g(z))>0$ for $z \in \mathbb{D}$. Class of close-to-convex functions denoted by $\mathcal{K}$. Functions in the class $\mathcal{K}$ are known to be univalent in $\mathbb{D}$. Geometrically, $f\in \mathcal{K}$
	means that the complement of the image-domain $f(\mathbb{D})$ is the union of rays that are
	disjoint (except that the origin of one ray may lie on another one of the rays).\\
	\vspace{3mm}

	In the late 1960's, Zalcman posed the conjecture that if  $f\in\mathcal{S}$ and is given by (\ref{p4_i001}), then
	\begin{equation}\label{p4_i005}
		|a_n^2-a_{2n-1}|\le (n-1)^2 \quad\mbox{ for } n\ge2,
	\end{equation}
	with equality  only for the Koebe function $k(z)=z/(1-z)^2$, or its rotation.
	It is important to note that the  Zalcman conjecture implies the celebrated Bieberbach conjecture $|a_n|\le n$ for $f\in\mathcal{S}$  
	(see \cite{Brown-Tsao-1986}), and a  well-known consequence of the area theorem  shows that  (\ref{p4_i005}) holds  for 
	$n=2$ (see \cite{Duren-book-1983}). The Zalcman conjecture remains an  open problem, even after de Branges' proof of the Bieberbach conjecture \cite{Branges-1985}.\\
	
	 For $f\in\mathcal{S}$, Krushkal  proved the Zalcman conjecture for $n=3$ (see \cite{Krushkal-1995}), and recently 
	for $n=4,5$ and $6$ (see \cite{Krushkal-2010}). For a simple and elegant proof of the Zalcman conjecture for the case $n=3$, see \cite{Krushkal-2010}. However the Zalcman conjecture for $f\in \mathcal{S}$ is still  open for $n>6$. On the other hand,  using complex geometry and universal Teichm\"{u}ller spaces, Krushkal   claimed  in an unpublished work \cite{Krushkal-Unpublished} to have proved the Zalcman conjecture  for all $n\ge 2$. Personal discussions with Prof. Krushkal indicates that there is a gap in the proof of  Krushkal's unpublished work \cite{Krushkal-Unpublished}, and so the Zalcman conjecture  remains open for the class $\mathcal{S}$ for $n>6$. In spite of these gaps, Krushkal's work \cite{Krushkal-Unpublished} contains  geometric and analytic ideas which may be  useful for  further study on the Zalcman conjecture and other similar  problems.\\

	 If $f\in \mathcal{S}$, then the coefficients $[f(z^2)]^{1/2}$ and $1/f(1/z)$ are polynomials in $a_n$, which contains the expression of the form $\lambda a_n^2-a_{2n-1}$,  pointed out by Pflunger \cite{Pfluger-1976}.\\
	
	The Zalcman conjecture and its other generalized form has been proved for some  subclasses of $\mathcal{S}$, such as starlike functions, typically real functions, close-to-convex functions \cite{Brown-Tsao-1986,Ma-1988,Ponnusamy-2017}. For basic properties of starlike functions, typically real functions and close-to-convex functions 
	we refer to \cite{Duren-book-1983}.\\

	In 1999 Ma \cite{Ma-1999} proposed a generalized Zalcman conjecture for $f\in \mathcal{S}$, namely that for $n\ge 2, m\ge 2$, $$|a_{n}a_{m}-a_{n+m-1}|\le (n-1)(m-1),$$
	which remains an open problem till date. However Ma  \cite{Ma-1999}  proved this generalized Zalcman conjecture for classes $\mathcal{S}^{*}$ and $\mathcal{S}_{\mathbb{R}}$, where   $\mathcal{S}_{\mathbb{R}}$  denotes 
	the class of 
	all functions in $\mathcal{S}$ with real coefficients.\\

The proofs of Theorems \ref{Thm 1}, \ref{Thm-2} are inspired by the work of Krushkal \cite{Krushkal-1995}. 
Roughly the idea is as follows:\\

Any function $f\in \mathcal{S}$ generates a family of holomorphic motion $f_{t}$ which is  a family of quasiconformal mappings. Also the coefficients of the functions in this family has the same homogeneity degree as the Zalcman functional as well as the generalized Zalcman functional. Since these maps are  a family of quasiconformal mappings, using the properties of extremal quasiconformal mappings we will see that for the extremal function of the generalized Zalcman functional the tangent map at $t=0$ is a $|\alpha|$-quasiconformal deformation of the Koebe's function, $|\alpha|\le 1$. More precisely, we will see that any extremal function can be written as the composition of a quasiconformal mapping $w^{\sigma}$ and an $|\alpha|$-quasiconformal deformation of the Koebe function. The goal is to show that $|\alpha|=1$, so that $w^{\sigma}$ reduced to a rotation. Thus we will deduce that the extremal function  is  a rotation of a Koebe function.\\

Consider the Koebe function
$$k(z)=\frac{z}{(1-z)^2}$$
and its rotations
$$k(z)=\frac{z}{(1-xz)^2},\,\,\,\, |x|=1$$		
and let 
$$k_{\alpha}(z)=\frac{z}{(1-\alpha xz)^2}, \ \mbox{ where }\ 0\le|\alpha|\le1.$$
	
	Before we state and prove our main results, we discuss some preliminary  ideas which will be useful in  our proofs.

	\section{Quasiconformal mappings}
	A homeomorphism $f$ on a domain $D\subset\mathbb{C}$ is called  K-quasiconformal if $f_{z}$ and $f_{\overline{z}}$, the partial derivatives with respect to $z$ and $\overline{z}$ in the sense of distribution, are locally in $\mathcal{L}^{2}$ and satisfy
	$$|\partial_{\overline{z}}f(z)|\le k|\partial_{z}f(z)|$$
	almost everywhere in $D$, where $k\in [0,1)$,  $k=(K-1)/(K+1)$, $f_{z}=\partial_{z}f$, $f_{\overline{z}}=\partial_{\overline{z}}f$, $\partial_{z}=\dfrac{1}{2}(\partial_{x}-i\partial_{y})$ and $\partial_{\overline{z}}=\dfrac{1}{2}(\partial_{x}+i\partial_{y})$.\\
	
The function $$\mu(z)=\frac{\partial_{\overline{z}}f}{\partial_{z}f} \,\,\,\,\,\,\,\,\,\,\,\mbox{ a.e. in } D$$
 is called the complex dilatation of the  map $f$. It is Borel-measurable (since the function $f$ is continuous) and satisfies the estimate 
	\begin{equation}
		|\mu(z)|\le \frac{K-1}{K+1}<1 \mbox{ and } \mbox{ ess sup }_{z\in D}|\mu(z)|=||\mu||_{\infty}\le k.
	\end{equation}


	{\bf The Beltrami equation:}
	 The equation 
	\begin{equation}
		\partial_{\overline{z}}f=\mu\,\partial_{z}f
	\end{equation}
	for a function $f$ in a domain $D$, where $\mu$ is a measurable bounded function in $D$ with  $\mbox{ ess sup}_{z\in D}\, |\mu(z)|=||\mu||_{\infty}< 1$ is called the Beltrami equation. The function $\mu$ in the equation is called the Beltrami differential. Thus $\mu$ belongs to the unit ball in the Banach space $\mathcal{L}^{\infty}(D)$. Consider
	
	$$K(f)=\frac{1+||\mu||_{\infty}}{1-||\mu||_{\infty}}$$ be the maximal dilatation of $f$. For the fundamental properties of quasiconformal mappings we refer to \cite{Ahlfors-1966},  \cite{Krzyz} and \cite{Lehto-1973}.\\
	
Consider two integral operators $S$ and $T$. The operator $S$ acts on functions $h\in \mathcal{L}^{p}$ for $p>2$ defined by
	
	\begin{equation}\label{Abhi-vasu-p3-eq1}
		Sh(z)=-\frac{1}{\pi} \int_{\mathbb{D}}\int h(\zeta)\left(\frac{1}{\zeta-z}-\frac{1}{\zeta}\right)d\xi d\eta, \,\,\,\,\, \zeta=\xi+i\eta.
	\end{equation}
	The operator $T$ is defined by
	$$Th(z)=-\frac{1}{\pi} \int_{\mathbb{D}}\int \frac{h(\zeta)}{(\zeta-z)^2}\,\,d\xi d\eta, \,\,\,\,\, \zeta=\xi+i\eta.$$
	
	  The operator $S$ is known as the Cauchy-Green integral operator, and the linear operator $T$ is known as the Hilbert transformation. The operator $T$ was initially defined only for function $\mathcal{C}_{0}^{2}$ ($\mathcal{C}_{0}$ with compact support). Functions of class $\mathcal{C}_{0}^{2}$ are dense in $\mathcal{L}^{2}$, so that using the isometry property, the operator $T$ extend to all of $\mathcal{L}^{2}$.  The operator $T$ is also correctly defined on functions $h\in \mathcal{L}^p$ with any $p>1$ and  bounded in the $\mathcal{L}^{p}$ norm
	$$||Th||_{\mathcal{L}^{p}}\le C_{p}||h||_{\mathcal{L}^p}.$$
	Moreover the bound $C_{p}$ can be chosen in such a way that $C_{p} \to 1$ for $p\to 2$.\\
	
	  It is well-known \cite[pg. 51]{Ahlfors-1966} that the operator $S$ satisfies 	
	$$|Sh(z)|\le K_p||h||_{p}|z|^{1-2/p}$$
	with a constant $K_{p}$ that depends only on $p$.\\
	\medskip

	For $h\in \mathcal{L}^p$, $p>2$, the relations  \cite[Lemma 3, pg. no. 90]{Ahlfors-1966},
	$$\partial_{\overline{z}}(Sh)=h$$
	$$\partial_{z}(Sh)=Th$$
	hold in the sense of distribution.

	\section{Holomorphic motions}
	  Holomorphic motions were introduced by Mane, Sad and Sullivan \cite{Sullivan-1983} in
	their study of the structural stability problem for  complex dynamical systems and  proved the first result in the topic which is called the $\lambda$-lemma. There is a strong connection between holomorphic motion and quasiconformal mappings, and the $\lambda$-lemma  illustrates this connection.
	After Mane, Sad and Sullivan had proved the $\lambda$-lemma, Sullivan and Thurston \cite{Sullivan-1986} proved an important result related to the extension of holomorphic motions. In the theory of holomorphic motions an important result is provided by  Slodkowski's lifting theorem, also known as the extended $\lambda$-lemma \cite{Slodkowski-1991}, which plays a vital role in the proof of our main result.\\
	
Basically a holomorphic motion is an isotopy of a subset $E$ of the extended complex
	plane $\mathbb{C}_{\infty}$ analytically parametrised by a complex variable $t$ in the unit
	disk $\mathbb{D}=\{z\in \mathbb{C}:|z|<1\}$. In this paper we use a special case of holomorphic isotopies generated by  analytic and univalent functions in the unit disk $\mathbb{D}$ which represents a special case of holomorphic motions.\\
	
	\begin{defn}
		 Let $E$ be a subset of $\mathbb{C}_{\infty}$ containing at least three points. A holomorphic motion of $E$ is a function $f:E \times \mathbb{D}\to \mathbb{C}_{\infty}$ such that 
		\begin{itemize}
			\item[(i)] for every fixed $z\in E$, the function $t \mapsto f(z,t):  E \times \mathbb{D}\to \mathbb{C}_{\infty}$ is holomorphic in $\mathbb{D}$;
			\item[(ii)] for every fixed $t\in \mathbb{D}$, the map $f(z,t)=f_{t}(z): E\to \mathbb{C}_{\infty}$ is injective;
			\item[(iii)] $f(z,0)=z$ for all $z\in E$.
		\end{itemize}
	\end{defn}
	
	The following  $\lambda$-lemma   assures that such a holomorphic motion  is a holomorphic family of quasiconformal maps.\\

	{\bf $\lambda$-lemma:} If $f:E\times \mathbb{D}\to\mathbb{C}_{\infty}$ is a holomorphic motion, then $f$ has an extension $\widetilde{f}:\overline{E}\times \mathbb{D}\to \mathbb{C}_{\infty}$ such that
	\begin{itemize}
		\item[(i)] $\widetilde{f}$ is a holomorphic motion of the closure $\overline{E}$ of $E$;
		\item[(ii)] each $\widetilde{f_{t}}(z)=\widetilde{f}(t,z):\overline{E}\to \mathbb{C}$ is quasiconformal on the interior of $\overline{E}$;
		\item[(iii)] $\widetilde{f}$ is jointly continuous in $(z,t)$.\\
	\end{itemize}
	 The obvious question  that arises  is  can we extend a holomorphic motion from any set to the whole sphere? The Slodkowski lifting theorem \cite{Slodkowski-1991} solves this problem which was posed by Sullivan and Thurston.\\
	
	{\bf Extended $\lambda$-lemma:} Any holomorphic motion $f: E\times \mathbb{D}\to \mathbb{C}_{\infty}$ can be extended to a holomorphic motion $\widetilde{f}:\mathbb{C}_{\infty}\times\mathbb{D}\to \mathbb{C}_{\infty}$, with $\widetilde{f}|_{E\times\mathbb{D}}=f.$\\
	
		 In view of \cite[Theorem 2]{Bers-1986}, the function $\phi:\mathbb{D}\to \mathcal{M}(\mathbb{D})$ defined by $$\phi(t)=\mu_{\widetilde{f_{t}}}(z)=\frac{\partial_{\overline{z}}\widetilde{f}(z,t)}{\partial_{z}\widetilde{f}(z,t)}$$ 
		 is holomorphic, where $\mathcal{M}(\mathbb{D})$ is a unit ball in $\mathcal{L}^{\infty}(\mathbb{D})$. 
		 
\section{Main results}
	
	 We  now state the main results of this paper. 
	
	\begin{thm}\label{Thm 1}
		
		Let $f\in \mathcal{S}$ and be given by (\ref{p4_i001}), then 
		\begin{equation}\label{eq-p4-1}
			|a_2a_3-a_4|\le 2
		\end{equation}
		with equality  only for  functions of the form
		$$\frac{z}{(1-e^{i \theta}z)^2},\,\,\, \mbox{ where }\,\, \theta \mbox{ is real}.$$
	\end{thm}
	\medskip
	
	\begin{thm}\label{Thm-2}
		
		Let $f\in \mathcal{S}$ and be given by (\ref{p4_i001}), then 
		\begin{equation}\label{eq-p4-2}
			|a_2a_4-a_5|\le 3 
		\end{equation}
		with equality  only for  functions of the form
		$$\frac{z}{(1-e^{i \theta}z)^2},\,\,\, \mbox{ where }\,\, \theta \mbox{ is real }.$$
	\end{thm}

	\section{Proof of Theorems \ref{Thm 1} and  \ref{Thm-2}}
The proofs of Theorems \ref{Thm 1}, \ref{Thm-2} are inspired by the work of Krushkal \cite{Krushkal-1995}. For  convenience we  convert the coefficients of functions in  $\mathcal{S}$ into the coefficients of functions in  $\sum$.\\
	
	 For $f\in \mathcal{S}$, we define a holomorphic motion $f^*$ of $\mathbb{D}\bigcup\{\infty \}$ {\it i.e.,}
	$f^* : \mathbb{D}\bigcup\{\infty \}\times \mathbb{D}\to \mathbb{C}_{\infty}$ by
	\[f^*(z,t)=\begin{cases}
		f_t(z)=\dfrac{1}{t}f(z,t)&\mbox{ if } z\in \mathbb{D}, \,\, t\in\mathbb{D}\\[5mm]
		\infty,&\mbox{ if } z=\infty,\,\, t\in \mathbb{D}.
	\end{cases}
	\]
	
	By the generalized $\lambda$-lemma, $f^*$ extends to a holomorphic motion say $\widetilde{f}$, $\widetilde{f}: \mathbb{C}_{\infty}\times\mathbb{D}\to\mathbb{C}_{\infty}$ such that $\widetilde{f}|_{(\mathbb{D}\bigcup\{\infty \})\times\mathbb{D}}=f^{*}$. {\it i.e.,} for fixed $t\in\mathbb{D}$, $\widetilde{f_{t}}(z)$ is quasiconformal on $\mathbb{C}_{\infty}$ and $\widetilde{f_{t}}|_{\overline{\mathbb{D}}}=f_{t}$. The Beltrami coefficient of $\widetilde{f_{t}}$ is given by
	
	$$\mu_{\widetilde{f_{t}}}(z,t)=\frac{\partial_{\overline{z}}\widetilde{f}(z,t)}{\partial_{z}\widetilde{f}(z,t)}.$$
	
		In view of \cite[Theorem 2]{Bers-1986}, since the mapping $t\mapsto \mu_{\widetilde{f_{t}}}$ is holomorphic, by Schwarz's lemma, we have
	$$||\mu_{\widetilde{f_{t}}}(z,t)||_{\infty}\le |t|.\\$$
	
	 Thus we have the following Taylor series expansion 
	
	\begin{equation}\label{TS}
		\mu_{\widetilde{f_{t}}}(z,t)=\mu_{1}(z)t+\mu_{2}(z)t^{2}+\mu_{3}(z)t^3+\cdots,
	\end{equation}

	where
	$$\mu_{\widetilde{f_{t}}}|_{\mathbb{D}}=0, \,\,\,\,\, \mbox{ since } \widetilde{f_{t}}|_{\mathbb{D}}=f_{t}.$$
	
	Similarly if $f\in \mathcal{S},\,\,\, \mbox{ then }$
	$$F(z)=\frac{1}{f\left(\dfrac{1}{z}\right)} \in \Sigma,\,\,\,\, F(\infty)=\infty.$$	
	
	 Therefore
	$$F\left(\dfrac{z}{t}\right)=\dfrac{1}{f\left(\dfrac{t}{z}\right)},\,\,\, F_{t}(z)=\dfrac{1}{f_{t}\left(\dfrac{1}{z}\right)}=\dfrac{t}{f\left(\dfrac{t}{z}\right)}=tF\left(\dfrac{z}{t}\right),$$
	
	and so
	\begin{equation}\label{b_n-1}
		F_{t}(z)=z+\sum_{n=0}^{\infty}b_{n}\frac{t^{n+1}}{z^n}.
	\end{equation}
	
	Now let $$g(z)=tF\left(\frac{z}{t}\right).$$
	Then we can define a holomorphic motion of $\mathbb{D^*}\cup \{0\}$ say  $F^*$, where $\mathbb{D^*}=\{z:|z|>1\}$, and   $F^*: \mathbb{D^*}\cup \{0\}\times\mathbb{D}\to\mathbb{C}_{\infty}$ is defined by 
	\[F^*(z,t)=\begin{cases}
		tF\left(\dfrac{z}{t}\right)&\mbox{ if } z\in \mathbb{D^*}, \,\, t\in\mathbb{D}\\[5mm]
		0,&\mbox{ if } z=0,\ \ t\in \mathbb{D}.
	\end{cases}
	\]
	 Similarly we can extend $F^*$ to a holomorphic motion of $\mathbb{C}_{\infty}$ say $\widetilde{F}$ by 
	$$\widetilde{F}(z,t)=\dfrac{1}{\widetilde{f}\left(\dfrac{1}{z},t\right)}$$
	with Beltrami coefficient 
	
	$$\mu_{\widetilde{F}(z,t)}=\mu_{\widetilde{f_{t}}(\frac{1}{z},t)}\frac{\overline{z}^2}{z^2}=\sum_{n=1}^{\infty}\mu_{n}\left(\frac{1}{z}\right)t^n\frac{\overline{z}^2}{z^2}=\widetilde{\mu_{1}}(z)t+\widetilde{\mu_{2}}(z)t^2+\widetilde{\mu_{3}}(z)t^3+\cdots.$$
	Thus $\widetilde{F}(z,t)$ satisfies
	$$\partial_{\overline{z}}\widetilde{F}(z,t)=\mu_{\widetilde{F}(z,t)}\partial_{z}\widetilde{F}(z,t),$$
	and we know that for each fixed $t$, $\widetilde{F}_{t}(z)$ is a quasiconformal map on $\mathbb{C}_{\infty}$ which is analytic on $\mathbb{D^*}$, and so by Weyl's lemma, $\partial_{\overline{z}}\widetilde{F}_{t}(z)=0$ for all $z\in \mathbb{D^*}$. Thus $\mu_{\widetilde{F}}$ has  compact support. Hence by the mapping theorem \cite[pg. 54]{Ahlfors-1966} these maps can be represented by $\widetilde{F_{t}}(z)=z+ S\partial_{\overline{z}}\widetilde{F_{t}}$, {\it i.e.},
	
	\begin{equation}\label{b_n}
		\widetilde{F_{t}}(z)=z-\frac{1}{\pi}\int_{\mathbb{D}} \int \partial_{\overline{z}}\widetilde{F_{t}}(\zeta)\left(\frac{1}{\zeta-z}-\frac{1}{\zeta}\right)\,d\xi d\eta,\,\,\, \mbox{ where }\, \zeta=\xi +\eta,
	\end{equation}
	and where 
	\begin{equation}
			\partial_{\overline{z}}\widetilde{F_{t}}=\mu_{\widetilde{F_{t}}}+\mu_{\widetilde{F_{t}}} T\mu_{\widetilde{F_{t}}} +\mu_{\widetilde{F_{t}}}  T\mu_{\widetilde{F_{t}}}  T\mu_{\widetilde{F_{t}}} +\mu_{\widetilde{F_{t}}}  T\mu_{\widetilde{F_{t}}}  T\mu_{\widetilde{F_{t}}}  T\mu_{\widetilde{F_{t}}} +\cdots.
	\end{equation}
	In the above equation, the right hand side series converges in $\mathcal{L}^p$, since the linear operator $h \mapsto T(\mu h)$ on $\mathcal{L}^p$ has norm less than or equal to $kC_{p}<1$. For $p>2$ sufficiently close to $2$.\\
	
	Now for each $n$, we will deduce the formula for $b_n$. In view of (\ref{b_n}), we have $|z|>1$, $|\zeta|<1$ and we have $|\zeta/z|<1$. From (\ref{b_n}) we have, 
	
	\begin{eqnarray*}
	\widetilde{F_{t}}(z)&=&z-\frac{1}{\pi}\int_{\mathbb{D}} \int \partial_{\overline{z}}\widetilde{F_{t}}(\zeta)\left(\frac{z}{\zeta(\zeta-z)}\right)\,d\xi d\eta\\ \nonumber
	&=&z+\frac{1}{\pi}\int_{\mathbb{D}} \int \partial_{\overline{z}}\widetilde{F_{t}}(\zeta)\left(\frac{1}{\zeta\left(1-\frac{\zeta}{z}\right)}\right)\,d\xi d\eta
	\end{eqnarray*}
which implies,
	\begin{eqnarray}\label{b_n-2}
		\widetilde{F_{t}}(z)&=&z+\frac{1}{\pi}\int_{\mathbb{D}} \int \partial_{\overline{z}}\widetilde{F_{t}}(\zeta)\frac{1}{\zeta}\left(1-\frac{\zeta}{z}\right)^{-1}\,d\xi d\eta\\ \nonumber
		&=&z+\sum_{n=0}^{\infty}\left(\frac{1}{\pi}\int_{\mathbb{D}} \int \partial_{\overline{z}}\widetilde{F_{t}}(\zeta)\zeta^{n-1}\,d\xi d\eta \right)\frac{1}{z^n}
	\end{eqnarray}
	By compairing (\ref{b_n-1}) and (\ref{b_n-2}), we obtain
	\begin{equation}\label{b_n formula}
		b_n=\frac{1}{\pi}\int_{\mathbb{D}} \int \partial_{\overline{z}}\widetilde{F_{t}}(\zeta)\zeta^{n-1}\,d\xi d\eta
	\end{equation}
Now write 
	\begin{eqnarray}\label{Beltrami}
		\partial_{\overline{z}}\widetilde{F_{t}}=\mu_{\widetilde{F_{t}}}\left(\partial_{z}\widetilde{F_{t}}-1\right) +\mu_{\widetilde{F_{t}}}.
	\end{eqnarray}
	Thus in view of equations (\ref{b_n formula}), (\ref{Beltrami}) and \cite[Lemma pg. no. 100]{Ahlfors-1966}, we have,
	
	\begin{equation}
		\widetilde{F_{t}}(z)=z-\frac{1}{\pi}\int_{\mathbb{D}} \int\mu_{\widetilde{F_{t}}}(\zeta)\left(\frac{1}{\zeta-z}-\frac{1}{\zeta}\right)\,d\xi d\eta+ O\left(||\mu_{\widetilde{F_{t}}}||_{\infty}^2\right),
	\end{equation}

	\begin{equation}
		b_n=\frac{1}{\pi}\int_{\mathbb{D}}\int \mu_{\widetilde{F_{t}}}(\zeta)\zeta^{n-1}\,d\xi d\eta+ O_n\left(||\mu_{\widetilde{F_{t}}}||_{\infty}^2\right).
	\end{equation} 
Where the remainder term is approximated uniformly on any compact subset of $\mathbb{C}$.\\

 We next mention a lemma (cf. \cite{Krushkal-1995} and \cite[Section 2.8]{Hayman-1987}) which we will apply to  the coefficients of  $\mu_{\widetilde{f_{t}}}$ for the fixed $z$.	
 \begin{lem}\label{lem} 	If $h$ is an analytic self-map of $\mathbb{D}$ fixing the origin with expansion $h(z)=\sum_{n=1}^{\infty} \alpha_nz^n$, then for $n\ge 1$,
 	$$|\alpha_{n+1}|\le 1-|\alpha_1|^2,$$
 \end{lem}
with equality only for $h(z)=z$.\\
\medskip

It therefore follows    from (\ref{TS}), and  Lemma \ref{lem}, that for a fixed $z$  
	\begin{equation}\label{inequality}
		|\mu_{n}(z)|\le 1- |\mu_{1}(z)|^2, \,\,\,\,\, n=2,3,4,\cdots
	\end{equation}
	
	for almost all $z\in \mathbb{D}$. 
	\medskip
	
	 We will need a much stronger version of  inequality (\ref{inequality}) as follows.

	\begin{equation}\label{star}
		||\mu_{n}||_\infty\le 1- ||\mu_1||_\infty^2.
	\end{equation}
	
	It is easy to see that  (\ref{star}) can be derived from (\ref{inequality}) provided 
	
	\begin{equation}\label{key}
		|\mu_{1}(z)|=\mbox{ constant a.e. on } \mathbb{D^*}.
	\end{equation}
	
	To establish (\ref{key}), Krushkal \cite{Krushkal-1995} proved a lemma known as the Surgery Lemma, which plays a vital role in the proof of our main result.

\begin{lem}[Surgery Lemma]
		Given a function $f\in \mathcal{S}$ with $a_{2}\ne0$, one can construct a holomorphic motion $f^{*}(z,t)$ on $\mathbb{C}_{\infty}\times\mathbb{D}$ such that for any $t$, the restriction of the fibre map $f_{t}^*(z)$ to $\mathbb{D}$ belongs to $\mathcal{S}$,
		and the family $f^*(z,t)$ satisfies the following properties: 
		\begin{itemize}
			\item[(i.)] $f^*(z,t)$ is also holomorphic in $t$, and its Beltrami coefficients
			\begin{equation}
				\mu_{\widetilde{f}_{t}^*(z)}=\mu_{1}^*(z)t+\mu_{2}^*(z)t^2+\mu_{3}^*(z)t^3+\cdots
			\end{equation}
		\end{itemize}
		satisfy
		\begin{equation}
			|\mu_{1}^*(z)|=constant \,\,\,\, \mbox{ on } \mathbb{D^*}
		\end{equation}
		\vspace{2mm}
		\begin{itemize}
			\item[(ii.)] $f_{t}^*(\infty)=\infty$ for all $t\in \mathbb{D}$
			\item[(iii.)] $a_{n}(f_{t}^*)=a_{n}(f_{t})+o_{n}(t^n), n=2,3,\cdots$.	
		\end{itemize}
			\end{lem}
		\vspace{3 mm}
 We  now give the proofs of Theorems \ref{Thm 1} and \ref{Thm-2}.

\subsection{Proof of Theorem \ref{Thm 1}} If $f\in \mathcal{S}$ any function maximizing
	 (\ref{eq-p4-1}), then $a_{2}(f)\ne 0$. If not, let $a_{2}=0$, then $b_{0}=0$ and so $ a_{2}a_{3}-a_{4}=b_{2}$. 
The  sharp inequality 
\begin{equation}\label{b_2}
|b_2|\le \frac{2}{3},
\end{equation}
was proved by Golusin \cite{Golusin-1938}. Jenkins \cite[Corollary 11]{Jenkins-1960}  gave another proof of this, and Duren \cite{Duren-1971} has given an elementary proof. Equality in (\ref{b_2}) holds only for the functions $$z\left(1+\frac{e^{i\theta}}{z^3}\right)^{\frac{2}{3}},\,\,\,\,\, \theta \mbox{ real},$$
which map $\{z:|z|>1\}$  onto the whole $w$-plane except for a cut consisting of three lines, each of length $2^{2/3}$, originating symmetrically from the origin. 
Thus we obtain the following sharp inequality
$$|b_2|\le \frac{2}{3},$$
 and so  it follows from (\ref{S to sigma}) that $|a_{2}a_{3}-a_{4}|=|b_{2}|\le 2/3$. This contradicts our hypothesis, since for  the Koebe function, $|a_{2}a_{3}-a_{4}|=2$, which is always greater than $2/3$  Hence  it is now clear that any function $f$ maximizing (\ref{eq-p4-1}) has second Taylor series coefficient $a_{2}\ne 0$, which leads us to use the Surgery Lemma. 
 \medskip
 
 So now we  consider the holomorphic motion 

	$$\widetilde{F}(z,t)=\frac{1}{\widetilde{f}(\frac{1}{z},t)}$$
		with Beltrami coefficient
		$$\mu_{\widetilde{F}(z,t)}=\mu_{\widetilde{f_{t}}\left(\frac{1}{z},t\right)}\frac{\overline{z}^2}{z^2}=\sum_{n=1}^{\infty}\mu_{n}\left(\frac{1}{z}\right)t^n\frac{\overline{z}^2}{z^2}=\widetilde{\mu_{1}}(z)t+\widetilde{\mu_{2}}(z)t^2+\widetilde{\mu_{3}}(z)t^3+\cdots$$
		generated by $\widetilde{f}$. From the proof of the Surgery Lemma (see \cite{Krushkal-1995}), we can assume that
		$$\widetilde{\mu_{1}}(z)=\alpha \widetilde{\mu}^0(z),$$
		where $$\widetilde{\mu}^0(z)=\frac{|z|}{\overline{z}},$$ 
		
and from \cite{Krushkal-1995} it is clear that for small $|t|$, we have the representation

$$\widetilde{F_t}(z)=w^{\sigma} \, o\, K_{\alpha t}(z)=z+b_0t+b_1t^2z^{-1}+\ldots,$$

where

	\[K_{\alpha t}(z)=\frac{1}{\tilde{k}_{\alpha t}\left(\frac{1}{z}\right)}=\begin{cases}
	z-2\alpha t+\alpha^2t^2z^{-1},& |z|>1,\\[5mm]
	z-2\alpha t|z|+\alpha^2t^2\overline{z},,& |z|<1,
\end{cases}
\]

has the Beltrami coefficient $-t\alpha\widetilde{\mu}^0$, $||\sigma(z,t)||_{\infty}\le |t|^2$, and the map

$$w^{\sigma}(\zeta)=\zeta+\sum_{n=0}^{\infty}b_n'\zeta^{-n}\,\,\,\, (b_n'=b_n'(t))$$
is conformal in the domain $K_{\alpha t}(\mathbb{D^*})\supset \{|\zeta|>1+O(t)\}$. 
\medskip

 We now  use the 
estimates of the coefficients of $w^{\sigma}$ mentioned in \cite{Krushkal-1995}, and obtain 

	\begin{equation}\label{Main}
		b_{0}'=O(t^2),\,\,\, b_{1}'=O(t^3),\,\,\, b_{2}'=O(t^4) \mbox{ and } |b_3'|\le \frac{(1-|\alpha|^2)^2}{\sqrt{3}(1-|t|)^2}|t|^4+O(t^6),
	\end{equation}
which gives
	$$b_{0}'b_{1}'=O(t^5).$$
	
	Since
	$$\widetilde{F_{t}}(z)=z-(2\alpha t-b_{0}')+\frac{\alpha^2t^2+b_{1}'}{z}+\frac{2\alpha tb_{1}'+b_{2}'}{z^2}+\frac{3\alpha^2t^2b_{1}'+4\alpha tb_{2}'+b_{3}'}{z^3}+\cdots,$$
 we have
	\begin{eqnarray}\label{Abhi-Vasu-00001}
		& b_{0}=& \frac{b_{0}'-2\alpha t}{t} \\[2 mm] \nonumber
		& b_{1}=& \frac{\alpha^2 t^2+b_{1}'}{t^2}\\[2 mm] \nonumber
		&b_{2}=&\frac{2\alpha tb_{1}'+b_{2}'}{t^3}\\[2 mm] \nonumber
		&b_{3}=&\frac{3\alpha^2t^2b_{1}'+4\alpha tb_{2}'+b_{3}'}{t^3},
	\end{eqnarray}
	
	and so from  (\ref{equ-p4-1}), we obtain
	$$	|a_{2}a_{3}-a_{4}|=|-b_{0}b_{1}+b_{2}|.$$

  Using  (\ref{Main}), after an easy computation, we obtain 
	
	$$	|a_{2}a_{3}-a_4|\le2|\alpha|^3, \,\,\,\,\,\,\,\, \mbox{ where } 0<\alpha\le 1,$$
	which implies that
	$$	|a_{2}a_{3}-a_4|\le 2,$$
	with equality  only if $|\alpha|=1$, {\it i.e.}, for the Koebe function and its rotations. This completes the proof of Theorem \ref{Thm 1}.
	
\vspace{3mm}
 \subsection{Proof of Theorem \ref{Thm-2}} If $f\in \mathcal{S}$ is a function maximizing (\ref{eq-p4-2}), then $a_{2}(f)\ne 0$. If not, let $a_{2}=0$, then $b_{0}=0$, and so $ a_{2}a_{4}-a_{5}=b_3-b_{1}^2$.  We need to find the sharp bound for $|b_{3}-b_1^2|$, and to do this we  use the method of Jenkins \cite{Jenkins-1960}.\\
 
  Jenkins \cite[Corollary 13]{Jenkins-1960} proved that, if $g\in \sum$, then for $\psi$ real and $0\le\sigma\le2,$

		\[\Re\left\{e^{-4i\psi}\left(b_3+\dfrac{1}{2}b_1^2-\sigma^2e^{2i\psi}b_1\right)\right\} \ge\begin{cases}
		-\dfrac{1}{2}-\dfrac{3}{16}\sigma^4+\dfrac{1}{8}\sigma^4\log\left(\dfrac{\sigma^2}{4}\right), \ 0<\sigma\le2\\[6 mm]
		-\dfrac{1}{2}, \,\,\, \sigma=0.
	\end{cases}
	\]
 Using the above inequality, we now find the sharp bound for $|b_3-b_1^2|$.\\

 Without loss in generality, choose $\psi$ real so that 
$$-\Re\left\{e^{-4i\psi}(b_3-b_1^2)\right\}=|b_3-b_1^2|,$$
and 
$$-\Re\left\{e^{-2i\psi}b_1\right\}\ge 0.$$
Then  it follows that
\begin{eqnarray*}
	|b_3-b_1^2|&\le& \frac{1}{2}+\frac{3}{16}\sigma^4-\frac{\sigma^4}{8}-\log\frac{\sigma^4}{4}+\frac{3}{2}\Re(e^{-4i\psi}b_1^2)-\sigma^2 \Re(e^{-2i\psi}b_1),\ \ 0<\sigma\le 2\\ \nonumber
	&\le& \frac{1}{2},\ \ \sigma=0.
\end{eqnarray*}
 It is easy to see that
 $$\Re\left\{e^{-4i\psi}b_1^2\right\}\le\left(\Re\left\{e^{-2i\psi}b_1\right\}\right)^2$$
with the strict inequality  unless
$$\Im\left\{e^{-2i\psi}b_1\right\}=0.$$
 Therefore
\begin{eqnarray*}
	|b_3-b_1^2|&\le& \frac{1}{2}+\frac{3}{16}\sigma^4-\frac{\sigma^4}{8}-\log\frac{\sigma^4}{4}+\frac{3}{2}\left(\Re\left\{e^{-2i\psi}b_1\right\}\right)^2-\sigma^2 \Re(e^{-2i\psi}b_1), \ \ 0<\sigma\le 2\\ \nonumber
	&\le& \frac{1}{2},\,\, \sigma=0.
\end{eqnarray*}

If $\Re\left\{e^{-2i\psi}b_1\right\}=0$, then taking $\sigma=0$ we obtain
$$|b_3-b_1^2|\le \frac{1}{2}.$$
If $\Re\left\{e^{-2i\psi}b_1\right\}>0$ then we can find $\sigma$ with $0<\sigma\le 2$, so that
$$\Re\left\{e^{-2i\psi}b_1\right\}=\frac{1}{4}\sigma^2\left(1-\log\left(\frac{\sigma^2}{4}\right)\right),$$
and using this we obtain
$$|b_3-b_1^2|\le \frac{1}{2}+\frac{1}{32}\sigma^4-\frac{1}{16}\sigma^4\log\frac{\sigma^2}{4}+\frac{3}{32}\sigma^2\log\frac{\sigma^2}{4}=\Theta(\sigma).$$
Clearly the function $\Theta(\sigma)$ is continuous on $0<\sigma\le2$, and further
$$|b_3-b_1^2|\le\max_{0<\sigma\le 2}\Theta(\sigma).$$
It is easy to see that the function $\Theta$ is increasing on $[0,\sigma]$  up to the point $\sigma=2e^{-\frac{1}{6}}$ and  is decreasing after  $\sigma=2e^{-\frac{1}{6}}$. Hence the maximum of $\Theta$ occurs at this value of $\sigma$ and so
 \begin{equation}\label{sharp}
 	|b_3-b_1^2|\le \frac{1}{2}+e^{-\frac{2}{3}}.
 \end{equation}
  Equality occurs in (\ref{sharp}) only for  functions of the form $h(z,2e^{-1/6},\psi)+k$, where $\psi$  is real, and $k$ is a constant, where $h(z,\sigma,\psi)$ is defined in \cite{Jenkins-1960}. Thus we have the sharp bound for $|b_3-b_1^2|$. 
 
  \medskip
  
  Thus 
  $$| a_{2}a_{4}-a_{5}|=|b_{3}-b_1^2|\le \frac{1}{2}+e^{-\frac{2}{3}},$$
  which is a contradiction, since for the Koebe function, $| a_{2}a_{4}-a_{5}|=3$. 
  \medskip
  
  Thus it is clear that any function $f$ maximizing (\ref{eq-p4-2}) has second Taylor series coefficient $a_{2}\ne 0$. So again we can use the Surgery Lemma and proceed as before.\\
 
Using equation (\ref{equ-p4-1}), we obtain
 
 $$| a_{2}a_{4}-a_{5}|=|-b_0b_2+2b_1b_0^2+b_{3}-b_{1}^2|.$$
Using equation (\ref{Main}) in above and after an easy computation, we obtain 

$$| a_{2}a_{4}-a_{5}|\le 3|\alpha|^4,\quad \mbox{ where } 0<\alpha\le 1.$$
This implies that
$$|\lambda a_{2}a_{4}-a_{5}|\le 3,$$
with equality  only if $|\alpha|=1$ {\it i.e.}, for the Koebe function and its rotations. This completes the proof of Theorem \ref{Thm-2}.

\vspace{4mm}
\noindent\textbf{Acknowledgement:} Authors thanks Professors T. Sugawa and D. K. Thomas  for their careful reading of this manuscript and giving  constructive suggestions for  improvements to this paper. The first author thanks SERB-CRG, and the second author thanks PMRF-MHRD, Govt. of India for their support.\\
	


\begin{thebibliography}{99}
		\bibitem{Ahlfors-1962}
	{\sc L. V. Ahlfors} and {\sc G. Weill G}, A uniqueness theorem for Beltrami equations, {\it Proc. Amer. Math. Soc.} {\bf 13} (1962), 975--978.
	
	
	\bibitem{Ahlfors-1966}
	{\sc L. V. Ahlfors}, Lectures on quasiconformal mappings. Van Nostrand, New York (1966).
	
		\bibitem{Bers-1986}
	{\sc Lipman Bers} and {\sc H. L. Royden}, Holomorphic families of injections, {\it Acta Math.} {\bf 157} (1986), 259--286.
	

	
	\bibitem{Branges-1985}
	{\sc L. de Branges}, A proof of the Bieberbach conjecture, {\it Acta Math.} {\bf 154} (1985), 137--152.
	
		\bibitem{Brown-Tsao-1986}
	{\sc J.E.Brown} and {\sc A.Tsao},  On the Zalcman conjecture for starlike and typically real functions,
	{\it Math. Z.} {\bf 191} (1986), 467--474.
	

	
	\bibitem{Conway-1995}
	{\sc J.B.Conway}, Functions of complex variables {\it II}, {\it Springer Verlag} (1995).
	
\bibitem{Duren-1971}
{\sc P.L.Duren}, Coefficients of meromorphic Schlicht functions, {\it Proc. Amer. Math. Soc.} {\bf 28} (1971), 169--172
	
	\bibitem{Duren-book-1983}
	{\sc P.L.Duren}, {\it Univalent functions} (Grundlehren der mathematischen Wissenschaften 259, New York, Berlin, Heidelberg, Tokyo),
	Springer-Verlag, 1983.
	
		\bibitem{Golusin-1938}
	{\sc G.M.Golusin}, Some evaluations of the coefficients of univalent functions. Mat. Sbornik {\bf 3} (1938), 321--330
	
	\bibitem{Hayman-1987}
	{\sc W.K.Hayman} and {\sc P. B. Kennedy}, Subharmonic functions, Vol. 1, {\it Academic Press, London}, 1987.
	
	\bibitem{Jenkins-1960}
{\sc J.A.Jenkins}, On certain coefficients of univalent functions. {\it Analytic functions. Princeton Univ. Press} (1960), 159--164.
	
	
	
	\bibitem{Krushkal-1995}{\sc S. L. Krushkal}, Univalent functions and holomorphic motions, {\it J. Analyse Math.} {\bf 66} (1995), 253--275.
	
	
	\bibitem{Krushkal-2010}
	{\sc S.L.Krushkal}, Proof of the Zalcman conjecture for initial coefficients, {\it Georgian Math. J.} {\bf 17} (2010), 663--681.
	
	\bibitem{Krushkal-Unpublished}
	{\sc S.L.Krushkal}, A short geometric proof of the Zalcman conjecture and Bieberbach conjecture, preprint, ariv:1408.1948.
	
	
		\bibitem{Krzyz}
	{\sc Julian Lawrynowicz} and {\sc Jan Krzyz}, Quasiconformal mappings in the Plane: Parametrical Methods, Springer, 1983.
	
	\bibitem{Lehto-1973}
	{\sc O. Lehto}, {\sc KI. Virtanen KI} Quasiconformal mapping in the plane. 2nd ed. New York and Heidelberg: Springer-Verlag; 1973.
	
	\bibitem{Lehto-1987}
	{\sc O. Lehto} Univalent functions and $Teichm\ddot{u}ller$ spaces. New york: Springer-Verlag; 1987.
	
		\bibitem{Ma-1988}
	{\sc W. Ma}, The Zalcman conjecture for close-to-convex functions, {\it Proc. Amer. Math. Soc.} {\bf 104} (1988), 741--744.
	
	\bibitem{Ma-1999}
	{\sc W. Ma}, Generalized Zalcman conjecture for starlike and typically real functions, {\it J. Math. Ana. Appl.} {\bf 234} (1999), 328--339.
	
	\bibitem{Ponnusamy-2017}
	{\sc Liulan Li} and {\sc S. Ponnusamy}, On the generalized Zalcman functional $\lambda a_n^2-a_{2n-1}$ in the close-to-convex family, {\it Proc. Amer. MAth. Soc.} {\bf 145} (2017), 833--846
	

	
	
	
	\bibitem{Pfluger-1976}
	{\sc Albert Pfluger}, On a coefficient problem for schlicht functions, {\it Advances in complex function
	theory (Proc. Sem., Univ. Maryland, College Park, Md., 1973), Lecture Notes in Math., Vol.
	505, Springer, Berlin}, 1976, 79-–91. 
	
	
	\bibitem{Sullivan-1986}
	{\sc D. Sullivan} and {\sc W. P. Thurston}, Extending holomorphic motions, {\it Acta Math.} {\bf 157} (1986), 243--257.
	

	
	
	
	
	\bibitem{Slodkowski-1991}
	{\sc Z. Slodkowski}, Holomorphic motions and polynomial hulls, {\it Proc Amer. Math. Soc.} {\bf 111}
	(1991), 347-355
	
	\bibitem{Sullivan-1983}
	{\sc R. Marl}, {\sc P. Sad} and {\sc D. Sullivan}, On dynamics of rational maps, {\it Ann. Sci. l~cole Norm.
	Sup.} {\bf 16} (1983), 193-217. 
	

	

	
\end{thebibliography}
\end{document}